\documentclass[12pt]{amsart}
\usepackage{amsfonts}

\theoremstyle{definition}

\theoremstyle{remark}

\newcommand{\const}{\mathop{\rm const}\limits}

\newcommand{\grad}{\mathop{\rm grad}\limits}

\begin{document}

\begin{center}

{\bf RIESZ'S AND BESSEL'S OPERATORS IN BILATERAL GRAND LEBESGUE SPACES} \par

\vspace{3mm}
{\bf E. Ostrovsky}\\

e - mail: galo@list.ru \\
 
\vspace{3mm}
{\bf E. Rogover}\\

e - mail: rogovee@gmail.com \\

\vspace{3mm}

{\bf L. Sirota}\\

e - mail: sirota@zahav.net.il \\

\vspace{3mm}
{\it Department of Mathematics and Statistics, Bar-Ilan University,
59200, Ramat Gan, Israel.}\\

\vspace{3mm}

 Abstract. \\
{\it In this paper we obtain the non - asymptotic estimations for Riesz's
and Bessel's potential integral operators in the so - called
Bilateral Grand Lebesgue Spaces. We also give examples to show the
sharpness of these inequalities.} \\

\end{center}

\vspace{3mm}


\vspace{3mm}


\vspace{3mm}

\section{Introduction}

\vspace{3mm}

 The linear integral operator $ I_{\alpha} f(x), $ or, wore precisely, the
 {\it family of operators} of a view

 $$
u(x) = I_{\alpha}f(x) = \int_{R^d} \frac{f(y) \ dy}{|x - y|^{d - \alpha}} \eqno(0)
 $$
is called {\it Riesz's integral operator,} or simply {\it Riesz's} potential, or
{\it fractional integral.} \par
Here  $ \alpha = \const \in (0,d) $  and $  |y| $ denotes usually  Euclidean norm
of the $ d - $ dimensional vector $ y; \ y \in R^d, \ d = 1,2,\ldots. $ \par
In the case $ d \ge 2 $ and $ \alpha = 1 \ u(x) $ coincides with the classical
Newton's potential.\par
 These operators are used in the theory of Fourier transform, theory of PDE,
 probability theory (study of potential functions for Markovian processes
  and spectral densities for stationary random fields), in the
 functional analysis, in particular, in the theory of interpolation of operators etc.,
 see for instance \cite{Bennet1}, \cite{Stein1}. \par

 We denote as usually the classical $ L_p $ Lebesgue norm

 $$
 |f|_p = \left( \int_{R^d} |f(x)|^p \ dx  \right)^{1/p};  \  f \in L_p \
 \Leftrightarrow |f|_p < \infty, \eqno(1)
 $$
and  denote $ L(a,b) = \cap_{p \in (a,b)} L_p.  $ \par
We will consider in the first two sections only the values $ p $ from the open interval
$ p \in (1,d/\alpha) $ and denote $ q = q(p) = pd/(d - \alpha p); $  or equally

$$
\frac{1}{q} = \frac{1}{p} - \frac{\alpha}{d}; \eqno(2)
$$
evidently, $ q \in (d/(d - \alpha), \infty). $ \par
The inverse function to the function $ q = q(p) $ has a view
$ p = p(q) = dq/(d + \alpha q );  $ note that if $ p \to 1 + 0 \ \Rightarrow
 q \to d/(d - \alpha) + 0 $ and if $ p \to d/\alpha - 0 \ \Rightarrow q \to \infty. $ \par

 More detail,

 $$
 p - 1 = \frac{(d - \alpha) (q - d/(d - \alpha))}{ d + \alpha \ q}
 $$
and
$$
\frac{d}{\alpha} - p = \frac{d^2}{\alpha(d + \alpha q) }.
$$

It is proved, see, e.g. in the book \cite{Stein1},  Chapter 5, p. 117 - 121 that
the following estimation holds for the Riesz's integral operator (0):

$$
|I_{\alpha} f|_{q} \le A(p)  \ |f|_p, \ \forall p \in (1,d/\alpha) \
 A(p) \in (0, \infty). \eqno(3)
 $$

 {\bf Our aim is a generalization of the estimation (3) on the so - called Bilateral
 Grand Lebesgue Spaces $ BGL = BGL(\psi) = G(\psi), $ i.e. when } $ f(\cdot)
 \in G(\psi). \ $ \par
  We recall briefly the definition and needed properties of these spaces.
  More details see in the works \cite{Fiorenza1}, \cite{Fiorenza2}, \cite{Ivaniec1},
   \cite{Ivaniec2}, \cite{Ostrovsky1}, \cite{Ostrovsky2}, \cite{Kozatchenko1},
  \cite{Jawerth1}, \cite{Karadzov1} etc. More about rearrangement invariant spaces
  see in the monographs \cite{Bennet1}, \cite{Krein1}. \par

\vspace{3mm}

For $a$ and $b$ constants, $1 \le a < b \le \infty,$ let $\psi =
\psi(p),$ $p \in (a,b),$ be a continuous positive
function such that there exists a limits (finite or not)
$ \psi(a + 0)$ and $\psi(b-0),$  with conditions $ \inf_{p \in (a,b)} > 0 $ and
 $\min\{\psi(a+0), \psi(b-0)\}> 0.$  We will denote the set of all these functions
 as $ \Psi(a,b). $ \par

The Bilateral Grand Lebesgue Space (in notation BGLS) $  G(\psi; a,b) =
 G(\psi) $ is the space of all measurable
functions $ \ f: R^d \to R \ $ endowed with  norm

$$
||f||G(\psi) \stackrel{def}{=}\sup_{p \in (a,b)}
\left[ \frac{ |f|_p}{\psi(p)} \right], \eqno(4)
$$
if it is finite.\par
 In the article \cite{Ostrovsky2} there are many examples of these spaces.
 For instance, in the case when  $ 1 \le a < b < \infty, \beta, \gamma \ge 0 $ and

 $$
 \psi(p) = \psi(a,b; \beta, \gamma; p) = (p - a)^{-\beta} (b - p)^{-\gamma};
 $$
we will denote
the correspondent $ G(\psi) $ space by  $ G(a,b; \beta, \gamma);  $ it
is not trivial, non - reflexive, non - separable
etc.  In the case $ b = \infty $ we need to take $ \gamma < 0 $ and define

$$
\psi(p) = \psi(a,b; \beta, \gamma; p) = (p - a)^{-\beta}, p \in (a, h);
$$

$$
\psi(p) = \psi(a,b; \beta, \gamma; p) = p^{- \gamma} = p^{- |\gamma|}, \ p \ge h,
$$
where the value $ h $ is the unique  solution of a continuity equation

$$
(h - a)^{- \beta} = h^{ - \gamma }
$$
in the set  $ h \in (a, \infty). $ \par

 The  $ G(\psi) $ spaces over some measurable space $ (X, F, \mu) $
with condition $ \mu(X) = 1 $  (probabilistic case)
appeared in \cite{Kozatchenko1}.\par
 The BGLS spaces
are rearrangement invariant spaces and moreover interpolation spaces
between the spaces $ L_1(R^d) $ and $ L_{\infty}(R^d) $ under real interpolation
method \cite{Carro1}, \cite{Jawerth1}. \par
It was proved also that in this case each $ G(\psi) $ space coincides
with the so - called {\it exponential Orlicz space,} up to norm equivalence. In others
quoted publications were investigated, for instance,
 their associate spaces, fundamental functions
$\phi(G(\psi; a,b);\delta),$ Fourier and singular operators,
conditions for convergence and compactness, reflexivity and
separability, martingales in these spaces, etc.\par

{\bf Remark 1.} If we introduce the {\it discontinuous} function

$$
\psi_r(p) = 1, \ p = r; \psi_r(p) = \infty, \ p \ne r, \ p,r \in (a,b)
$$
and define formally  $ C/\infty = 0, \ C = \const \in R^1, $ then  the norm
in the space $ G(\psi_r) $ coincides with the $ L_r $ norm:

$$
||f||G(\psi_r) = |f|_r.
$$

Thus, the Bilateral Grand Lebesgue spaces are direct generalization of the
classical exponential Orlicz's spaces and Lebesgue spaces $ L_r. $ \par

The BGLS norm estimates, in particular, Orlicz norm estimates for
measurable functions, e.g., for random variables are used in PDE
\cite{Fiorenza1}, \cite{Ivaniec1}, theory of probability in Banach spaces
\cite{Ledoux1}, \cite{Kozatchenko1},
\cite{Ostrovsky1}, in the modern non-parametrical statistics, for
example, in the so-called regression problem \cite{Ostrovsky1}.\par

The article is organized as follows. In the next section we obtain
the main result: upper bounds for Riesz's operators in the Bilateral
Grand Lebesgue spaces. In the third section we construct some examples
in order to illustrate the precision of upper estimations. \par

 In the fourth section we  investigate the so - called Bessel's potential
 operator and establish its boundedness in the Grand Lebesgue Spaces.\par

  In the fifth section we  investigate the so - called truncated Riesz's
 operator and establish its boundedness in the Grand Lebesgue Spaces,
study the sharpness of the
obtained results by the building of the suitable examples.\par
 The last section contains some slight generalizations of obtained results.\par

\vspace{3mm}

 We use symbols $C(X,Y),$ $C(p,q;\psi),$ etc., to denote positive
constants along with parameters they depend on, or at least
dependence on which is essential in our study. To distinguish
between two different constants depending on the same parameters
we will additionally enumerate them, like $C_1(X,Y)$ and
$C_2(X,Y).$ The relation $ g(\cdot) \asymp h(\cdot), \ p \in (A,B), $
where $ g = g(p), \ h = h(p), \ g,h: (A,B) \to R_+, $
denotes as usually

$$
0< \inf_{p\in (A,B)} h(p)/g(p) \le \sup_{p \in(A,B)}h(p)/g(p)<\infty.
$$
The symbol $ \sim $ will denote usual equivalence in the limit
sense.\par
We will denote as ordinary the indicator function
$$
I(x \in A) = 1, x \in A, \ I(x \in A) = 0, x \notin A;
$$
here $ A $ is a measurable set.\par
 All the passing to the limit in this article may be grounded by means
 of Lebesgue dominated convergence theorem.\par

\bigskip

\section{Main result: upper estimations for Riesz potential}

\vspace{3mm}

{\bf A.} Estimations for Riesz potential. \par

\vspace{3mm}

Let $ \psi(\cdot) \in \Psi(a,b), \ $ where $ 1 \le a < b \le d/\alpha. $  Define for the
arbitrary function  $ \psi \in \Psi(a,b)  $ the auxiliary function of the {\it variable}
$ q, \ q \in (d/(d - \alpha), \infty): $

$$
\psi_{\alpha, d}(p) = [(p-1) \ (d/\alpha - p)]^{-1 + \alpha/d } \ \psi(p), \eqno(5)
$$
where $ p $ is the following function of variable $ q:$

$$
p = p(q) =  \frac{d q}{d + \alpha q },
$$
or up to equivalence:

$$
\psi_{\alpha, d}(p(q)) \asymp \zeta_{\alpha,d}(q) = \zeta_{\alpha,d, \psi(\cdot) }(q)
 \stackrel{def}{=} \left[ \frac{q^2 }{ q - d/(d - \alpha) } \right]^{1 - \alpha/d} \
\psi \left( \frac{d q}{d + \alpha q }  \right). \eqno(6)
$$
Notice that the values $ p = 1 + 0 $ and $ p = d/\alpha - 0 $ are critical points;
another points are not interest. \par

\vspace{3mm}

{\bf Theorem 1.} Let $ f \in G(\psi). $ Then

\vspace{3mm}

$$
||I_{\alpha} \ f||G(\zeta_{\alpha, d}) \le C( \alpha, d, \psi) \
|| f ||G(\psi).  \eqno(7).
$$

\vspace{3mm}

{\bf Example 1.} When $ a = 1, b = d/\alpha,  \beta, \gamma > 0, $ and
$  f \in G(a,b; \beta, \gamma), \ f \ne 0, $ then

$$
I_{\alpha} f(\cdot) \in G( d/(d - \alpha), \infty;
 \beta + 1 - \alpha/d), -\gamma - 1 + \alpha/d ). \eqno(8)
$$

\vspace{3mm}

 {\bf Lemma 1.} {\it The constant $  A(p) $ from the inequality (3) may be estimates
 as follows:}

$$
A(p) \le C_1(d,\alpha) \frac{ \alpha^{-1} (d - \alpha)^{-1} }
{ [(p - 1) \ (d/\alpha - p)]^{1 - \alpha/d} }, \ p \in(1, d/\alpha),  \eqno(9)
$$

where $ C_1(d, \alpha) $ is bounded continuous function of the variable $ \alpha $
on the closed interval $ \alpha \in [0, d].$ \par
{\bf Proof } of lemma 1. \par
 It follows from the mentioned monograph of E.M.Stein
\cite{Stein1}, Chapter 5, pp. 117 - 121 and from the \cite{Bennet1}, Chapter 4,
section 4, pp. 216 - 230  a weak proposition,  after some calculations. \par

 In detail, let us denote

 $$
 \phi(x) = |x|^{\alpha - d};
 $$
then the Riesz's  transform may be rewritten as a convolution
$$
I_{\alpha}f = f*\phi. \eqno(10)
$$

 The distribution function $m_{\phi}(\lambda) $ for the positive
  function $ \phi $ has a view

 $$
 m_{\phi}(\lambda) = C_2(d, \alpha) \lambda^{d/(d - \alpha)}, \ \lambda \in (0,\infty).
 \eqno(11)
 $$
 Therefore, the function $ \phi(\cdot) $ belongs to the generalized Lebesgue space
 $ L_{d/(d - \alpha), \infty} $ and following the operator $ I_{\alpha} $ is of weak
 type $ (1, d/(d - \alpha)):  $

$$
|I_{\alpha} f|_{d/(d - \alpha), \infty } \le |\phi|_{d/(d - \alpha)} \ |f|_1
\stackrel{def}{=} M_0 \ |f|_1. \eqno(12)
$$
 The operator $ I_{\alpha} $ is also the the weak type $ (d/(\alpha), \infty): $

$$
|I_{\alpha} f| \le M_1 \ \int_0^{\infty}f^*(t) t^{\alpha/d - 1} dt,  \eqno(13)
$$
where $ f^* $ denotes the inversion function to the distribution function for
the function $ f. $  More exact calculation show us that

$$
\max(M_0, M_1) \le C(d)/ \min(\alpha, d - \alpha).
$$

 We will use the classical Marcinkiewicz's interpolation theorem, see e.g.
 \cite{Bennet1}, Chapter 4, corollary 4.14, p. 225.   Finding the
 parameter $ \theta $ from the equation

 $$
 \frac{1}{p} = \frac{1 - \theta}{1} + \frac{\theta}{d/\alpha},
 $$
we have:
$$
\theta = \frac{d }{d - \alpha} \frac{p-1}{p}
$$
and correspondingly
$$
1 - \theta = \frac{\alpha}{d - \alpha} \frac{d/\alpha - p}{p}.
$$
From the Marcinkiewicz's interpolation theorem follows that the operator $ I_{\alpha} $
is bounded as operator from the space $ L_p $ into the space $ L_q $ with the  following norm estimation:

$$
||I_{\alpha}||(L_p \to L_q) \le C_2(d,\alpha) \frac{\max(M_0,\ M_1)}{\theta(1 - \theta)}
\le C_3(d,\alpha) (p-1)^{-1} \ (d/\alpha - p)^{-1}.
$$
 But in the book \cite{Adams1}, pp. 49 - 54 there is described
 other approach, which used the so - called maximal operator

 $$
 M f(x) = \sup_{r > 0} \left[ r^{-d} \ \int_{y: |y - x| \le r} |f(y)| \ dy \right].
 $$
It is known ( \cite{Stein1}, p. 173 - 188) that

$$
| M f|_p \le C(d) \frac{p}{p-1} \ |f|_p, \ p \in (1, \infty),
$$
and if  the variable $ p $ changed in the finite interval $ p \in (1, d/\alpha), $
then

$$
| M f|_p \le C(\alpha, d) \frac{1}{p-1} \ |f|_p.
$$
 The exact computation following \cite{Adams1}  show us that

 $$
 |I_{\alpha} f|_q \le \frac{C \ |f|_p }{[(p - 1) \ (d/\alpha - p)]^{1 - \alpha/d} }.
 \eqno(14)
 $$

 Further we will show that the last estimation (14) is sharp  and will
 prove more general proposition. \par

\vspace{3mm}

{\bf Proof} of the theorem 1. Denote for the simplicity $ u = I_{\alpha}f; \
u: R^d \to R. $\par
 We can assume without loss of generality that $ ||f|| G(\psi) = 1; $
this means that

$$
\forall p \in (a,b) \ \Rightarrow |f|_p \le \psi(p).
$$

Using the inequality (3) with the estimation (14), we obtain:

$$
|u|_q \le A(p) \ \psi(p) \le C_4(d,\alpha) \
[(p-1) \ (d/\alpha - p)]^{-1 + \alpha/d} \ \psi(p)
$$

$$
\le C_4(d,\alpha) \ \psi_{\alpha, d}(p) \ ||f||G(\psi). \eqno(15)
$$
The assertion of theorem 1 follows after replacing $ p = dq/(d + \alpha q), $
 dividing on the  $\zeta_{\alpha, d}(q) $ and on the basis of the definition of the
$ G(\psi) $ spaces.

\vspace{3mm}

{\bf B.} Derivatives of Riesz's potential.\par

\vspace{3mm}

Let $ \xi = \vec{\xi} = (\xi_1, \xi_2, \ldots, \xi_d) $ be a non - negative integer
multiindex: $ \xi_j \ge 0, \ \xi_j = 0,1,2, \ldots. $  We define for such a indices
$ |\xi| = \sum_{j=1}^d \xi_j. $ \par
We restrict the values $ \xi $ as follows: $ | \xi| < \alpha $ and denote

 $$
 \alpha(\xi) = \alpha - |\xi|; \ \alpha(0) = \alpha; \ \alpha(\xi) \in (0, \alpha];
 $$

$$
u^{(\xi)}(x) \stackrel{def}{=}
D^{(\xi)} u(x) =  \frac{\partial^{\xi} u }
{ \partial x_1^{\xi_1} \ \partial x_2^{\xi_2} \ldots \partial x_d^{\xi_d} },
$$
here as before $ u(x) = I_{\alpha}f(x). $ \par
 In this subsection we assume that the value $ p $ belongs to the interval
 $ p \in (1, d/\alpha(\xi) ) = (1, d/(\alpha - |\xi|)) $ and denote

 $$
 q^{(\xi)} = q^{(\xi)}(p) = \frac{d p}{ d - (\alpha - |\xi|)p } \in
 \left( \frac{d }{d - \alpha(\xi)}, \ \infty \right);
 $$

$$
A^{(\xi)}(p) =  \frac{ \alpha^{-1} (d - \alpha)^{-1} }
{ [(p - 1) \ (d/\alpha(\xi) - p)]^{1 - \alpha(\xi)/d} }, \ p \in(1, d/\alpha(\xi)).
\eqno(16)
$$

\vspace{3mm}

{\bf Theorem 2.}

\vspace{3mm}

$$
\left| \ u^{(\xi)} \ \right|_q \le C \ A^{(\xi)}(p) \ |f|_p, \eqno(17)
$$
and if $ f \in G(\psi(1, d/\alpha(\xi))), $ then
$ \left| \ u ^{(\xi)} \right|_{q} \le $

$$
C \ ||f||G\Psi \ \times
\left[ \frac{q}{d/(d - \alpha + |\xi|) } \right]^{1 - (\alpha - |\xi|)/d } \times
\psi \left( \frac{dq}{d + q(\alpha - |\xi|)}  \right). \eqno(18)
$$
{\bf Proof.} We follow the \cite{Adams1}, p. 58 - 59.  The expression for the
function $ u = I_{\alpha}f $ may be rewritten as the convolution
$ u = I_{\alpha}*f, $ where $ I_{\alpha}(x) = |x|^{d - \alpha}. $ Therefore,

$$
D^{(\xi)} (I_{\alpha}*f) = (D^{(\xi)}I_{\alpha} * f).
$$
 The estimation

 $$
  \left| D^{(\xi)} \ I_{\alpha}(x)  \right| \le A(\xi, d, \alpha) \
  I_{\alpha - |\xi|}(x) \eqno(19)
 $$
see in \cite{Adams1}, p. 57; hence

$$
 \left|  u^{(\xi)} \right| \le C \ I_{\alpha - |\xi|}*f. \eqno(20)
$$
 The end of the proof of theorem 2 is at the same as the proof of theorem 1.\par
{\bf Example 2.} Newton potential. \par
Recall that if $ d \ge 3 $ and  $  f(\cdot) \in \cup_{p \in (1,d/2)}L_p(R^d), $
then

$$
f(x) = \frac{-1}{(d-2) \ \omega(d-1)} \int_{R^d} \frac{\Delta f(y) \ dy }
{|x-y|^{d-2} };
$$
here $ \Delta $ denotes the Laplace operator and
$ \omega(d-1) = 2 \pi^{d/2}/ \Gamma(d/2), $
or

$$
f(x) = \frac{1}{\omega(d-1)} \int_{R^d} \frac{(\grad f(y) \cdot (x - y)) \ dy}
{|x - y|^d }.
$$
We obtain noticing that

$$
f(x) = C(d) \ I_{2}( \Delta f), \ f(x) = C_2 \ I_1( ( \grad f(x) \cdot x) ) ):
$$

$$
|f|_q \le C_1(d) \frac{|\Delta \ f|_p }{ \left[(p-1)(d/2 - p) \right]^{1 - 2/d} }, \
q = 1/(1/p - 2/d), \ p \in (1, d/2).
$$
or

$$
|f|_q \le C_2(d) \frac{| \grad \ f|_p }{\left[(p-1)(d - p) \right]^{1 - 1/d} }, \
q = 1/(1/p - 1/d), \ p \in (1,d).
$$

\hfill $\Box$

\bigskip

\section{Low bounds for Riesz potential.}

\vspace{3mm}

In this section we built some examples in order to illustrate the
 exactness of upper
estimations. We consider only the case $ a = 1, \ b = d/\alpha; $  another cases
are trivial.\par
 It is sufficient for constructing  the low estimations
 to consider only the one - dimensional case: $ d = 1,$
i.e. $ x,y \in R^1. $ \par
 Let us denote for the mentioned  values $ p, q(p), f \in L_p, f \ne 0 $
$$
V_{\alpha,d}(f,p) = V(f,p) = \frac{|I_{\alpha} f |_q \cdot
[(p-1) \ (d/\alpha - p)]^{1 - \alpha/d} }{|f|_p }, \eqno(21)
$$
where as before $ p \in (1,d/\alpha), \ q = q(p) = pd/(d - \alpha p) \in
(d/(d - \alpha),\infty). $ \par
From the inequality (3) with  the concrete values of parameters from (14)
 follows that for {\it some } non - trivial values $ C^{(1)} = C^{(1)}(\alpha,d),
 C^{(2)} = C^{(2)}(\alpha,d) $

$$
 \sup_{f \in L(1,d/\alpha), f \ne 0} \overline{\lim}_{p \to 1+0}
 V(f,p) \le C^{(1)},
$$

$$
 \sup_{f \in L(1,d/\alpha), f \ne 0} \overline{\lim}_{p \to d/\alpha - 0}
 V( f,p) \le C^{(2)}.
$$

We intend to prove an inverse inequality at both the critical points $ p \to
 1 + 0 $ and $ p \to d/\alpha - 0. $ \par

\vspace{3mm}

{\bf Theorem 3.}  For all the values $ \alpha \in (0,1) $
there exist a constants $ C_1(\alpha), C_2(\alpha) \in (0,\infty) $ for which

\vspace{3mm}

$$
 \sup_{f \in L(1,1/\alpha), f \ne 0} \underline{\lim}_{p \to 1+0}
 V_{\alpha,1}(f,p) \ge C_1, \eqno(22)
$$

$$
 \sup_{f \in L(1,1/\alpha), f \ne 0} \underline{\lim}_{p \to 1/\alpha - 0}
 V_{\alpha,1}(f,p) \ge C_2. \eqno(23)
$$

\vspace{3mm}

{\bf Proof.} First of all we consider the case $ p \to 1 + 0; $ for definiteness
we restrict the value p in the semi - closed interval
$ p \in (1, 0.5( 1 + 1/\alpha)]. $ \par
  Let us consider a function (more exactly, a family of the functions) of a view
$$
g(x) = g_{\Delta}(x) = x^{-1} |\log x|^{\Delta}, \ x \in (e, \infty)
$$
and $ g_{\Delta}(x) = 0 $  otherwise.  Here the value $ \Delta $ is fixed
constant number $ \ \Delta = \const \ge 0.$ \par
 We have:

$$
|g_{\Delta}|_p^p = \int_{e}^{\infty} x^{-p} \ (\log x)^{\Delta p} \ dx =
\int_1^{\infty} \exp (-(p-1)y) \ y^{\Delta p} \ dy \sim
$$

$$
\int_0^{\infty} \exp (-(p-1)y) \ y^{\Delta p} \ dy
= (p-1)^{- \Delta p - 1} \ \Gamma(\Delta p + 1), \eqno(24)
$$
where $ \Gamma(\cdot) $ denotes usually Gamma - function. \par
 It follows from the equality (24) that as $ p \to 1 + 0 $

$$
 |g_{\Delta}|_p  \asymp \ (p - 1)^{-\Delta - 1} \asymp
 [1/(1 - \alpha)- q]^{\Delta + 1}.
$$
Further, let us investigate the behavior of the corresponding function
$ u_{\Delta} = I_{\alpha} g_{\Delta} $ as $ x  \to \infty, \ x > 1: $

$$
u_{\Delta}(x) = \int_e^{\infty} \frac{y^{-1} \ (\log y)^{\Delta} \ dy }
{|x - y|^{1 - \alpha}} =
x^{\alpha - 1} \int_{e/x}^{\infty} \frac{z^{-1} \ (\log x + \log z)^{\Delta} \ dz }
{ |z - 1|^{1 - \alpha}} \sim
$$

$$
x^{\alpha - 1} \ (\log x)^{\Delta} \int_{1/x}^{\infty} z^{-1} \
|z - 1|^{\alpha - 1} \ dz \sim
x^{\alpha - 1} \ (\log x)^{\Delta} \int_{1/x}^{1/e} z^{-1} \
|z - 1|^{\alpha - 1} \ dz \sim
$$

$$
x^{\alpha - 1} \ (\log x)^{\Delta} \int_{1/x}^{1/e} z^{-1} \ dz \sim
x^{\alpha - 1} \ (\log x)^{\Delta + 1};
$$
therefore as $ q \to 1/(1 - \alpha) + 0 \  $

$$
|u_{\Delta}|_q^q \asymp
 \int_e^{\infty} x^{q(\alpha - 1)} \ (\log x)^{q(\Delta + 1)} \ dx \sim
 \int_1^{\infty} x^{q(\alpha - 1)} \ (\log x)^{q(\Delta + 1)} \ dx =
$$

$$
\left[q - \frac{1}{1 - \alpha} \right]^{-q(\Delta + 1) - 1} \ \Gamma(q(\Delta+1) + 1),
$$
 following

 $$
 |u_{\Delta}|_q  \asymp \left[q - \frac{1}{1 - \alpha} \right]^{-\Delta -1  - 1/q} \asymp
 \left[q - \frac{1}{1 - \alpha} \right]^{-\Delta -2 + \alpha} \asymp
 (p-1)^{-\Delta - 2 + \alpha}.
 $$

 Substituting into the expression for the functional $ V_{\alpha,1}, $ we obtain for
 some constant $ C_1 = C_1(\alpha,1, \Delta) $ and for the value $ p $ tending  to
 $ 1 + 0: $

 $$
 V_{\alpha,1}( g_{\Delta},p) \asymp
 \frac{(p-1)^{-\Delta - 2 + \alpha} \ (p-1)^{1 - \alpha} }{ (p-1)^{1 - \alpha} }
\asymp C_1.
 $$

We conclude for all the the values $ \Delta \ge 0 $

$$
 \underline{\lim}_{p \to 1 + 0} V_{\alpha, 1} (g_{\Delta},p) \ge C_1.
$$

 We prove now the second assertion of theorem 3.
  Let us consider now the case $ p \to 1/\alpha - 0,  $ then $ q \to \infty.  $
We introduce as an example the family of a functions

$$
f_{\Delta}(x) = x^{-\alpha} \ |\log x|^{\Delta}, \ x \in (0,1/e)
$$
and $ f_{\Delta}(x) = 0 $ otherwise. As before, $ \Delta $ is arbitrary non - negative
constant parameter.\par

We  calculate:

$$
|f_{\Delta}|_p^p = \int_0^{1/e} x^{-\alpha p} \ |\log x|^{\Delta p} \ dx \asymp
C(\alpha) \ (1/\alpha - p)^{- \Delta p - 1};
$$

$$
|f_{\Delta}|_p \asymp  (1/\alpha - p)^{- \Delta - \alpha}.
$$
 Further,  let us denote $ v_{\Delta}(x) = I_{\alpha} f_{\Delta}(x).$  We have as
 $ x \to 0+, x \in (0, 1/e):$

$$
v_{\Delta}(x) = \int_0^{1/e}
\frac{y^{-\alpha} \ |\log y|^{\Delta} }{ |x - y|^{1 - \alpha} }
\ dy \sim |\log x|^{\Delta} \int_0^{1/ex} \frac{z^{-\alpha} \ dz}{|z - 1|^{1 - \alpha }}
$$

$$
\asymp  |\log x|^{\Delta} \int_e^{1/ex} \frac{z^{-\alpha} \ dz}{|z - 1|^{1 - \alpha }} \asymp |\log x|^{\Delta} \int_e^{1/ex} dz/z \asymp |\log x|^{\Delta + 1};
$$

$$
\left|v_{\Delta} \right|_q^q \asymp \int_0^1 |\log x|^{q(\Delta + 1) } \ dx =
\Gamma(q(\Delta + 1) + 1),
$$

$$
\left|v_{\Delta} \right|_q  \asymp q^{\Delta + 1}.
$$

 We  used the Stirling's formula as $ q \to \infty $ and correspondingly
  $ p \to 1/\alpha - 0. $ \par

 Substituting into the expression for the functional $ V_{\alpha,1}, $ we obtain for
 suitable constant $ C_2 = C_2(\alpha,1, \Delta) $ and for the value $ p $  tending to
 $ 1/\alpha - 0: $

 $$
 V_{\alpha,1}(f_{\Delta},p) \asymp \frac{(1/\alpha - p)^{-\Delta - 1} \
 (1/\alpha - p)^{1 - \alpha} }{(1/\alpha - p)^{-\Delta - \alpha} }
  \asymp C_2.
 $$

This completes the proof of the second assertion of theorem 3. \par
\hfill $\Box$

\vspace{3mm}

{\bf Remark 2.} If we want to construct a single example for the proof
of theorem 3, we can choose the function

$$
h_{\Delta}(x) = f_{\Delta}(x) + g_{\Delta}(x).
$$
{\bf Corollary 1.} If we choose as a function $ \psi(p) = \psi_{\Delta}(p) $ as
follows:

$$
\psi_{\Delta}(p) = |h_{\Delta}|_p, \ p \in (1, d/\alpha),
$$
we see that as long as the supports of a functions $ f_{\Delta} $ and $ g_{\Delta} $
are disjoint, for both the sets of a values $ p: \ p \in (1,0.5(1 + d/\alpha)) $
and $ p \in (0.5(1 + d/\alpha), d/\alpha)$

$$
V_{\alpha,1}(h_{\Delta}, p) \asymp C.
$$
Therefore, the proposition of theorem 1 is unprovable. \par

\vspace{3mm}

{\bf Remark 3.} The "extremal" case $ \alpha = 0, $ in which the integral operator
$ R_j f(x) $ understood as a principal value with {\it "signed kernel"}:

 $$
 R_j f(x) = p.v. \int_{R^d} \frac{x_j - y_j }{|x - y|^{d + 1} } \ dy,
 $$
  and is called a Riesz's
singular integral transform, is considered in \cite{Ostrovsky2};
there is obtained the exact $ G(\psi) $ estimations for  the norm of $ R_j $ and is
showed its exactness.\par
 Namely, in the book \cite{Stein1}, pp. 26 - 48 and 51 is proved the following estimation:

 $$
 |R_j f|_p \le K(p) \ |f|_p, \ p \in (1, \infty), \eqno(25)
 $$
where

$$
K(p) \le C_5(d) \frac{p^2}{ p - 1}. \eqno(26)
$$
 Notice that the inequalities (25) and (26) are "extremal"  cases for the estimations
 for the Riesz's transform: $ \lim_{\alpha \to 0+} q(p) = p $ and the interval for the
 values $ p $ tends to whole semi - axis $ (1, \infty). $ \par
The estimations (25) and (26) are true for more general singular integral operators
and more generally for pseudodifferential operators.\par
 From the inequalities (25) and (26) follows that if $ f \in G(\psi; a,b),  $ where
$ 1 \le a < b \le \infty, $ then $ R_j f \in G(\psi^{(1)}), $  where

$$
\psi^{(1)}(p) = \frac{p^2}{ p - 1} \cdot \psi(p):
$$

$$
||R_j f||G(\psi^{(1)}) \le C_5(d) \ ||f||G(\psi), \eqno(27)
$$
and the last estimation (27) is exact still in the one - dimensional case,
see \cite{Ostrovsky2},  in which
the Riesz's transform coincides with the classical Hilbert transform. \par

\bigskip

\section{The case of Bessel's potential}

\vspace{3mm}

{\bf A.} Estimation of Bessel's potential.

\vspace{3mm}

We consider here a so - called Bessel's integral operator or equally Bessel potential:
$$
v(x) = L_{\alpha}f = G_{\alpha}*f(x), \ G_{\alpha}(x) = |x|^{-(d - \alpha)/2}
K_{(d - \alpha)/2}(|x|),
$$
where $ K_{\xi}(x) $ is a modified Bessel's function of the third kind with
index $ \xi, \ \xi \ge 0, $  or a Macdonald's  function. \par
 It is true the following representation for the function $ G_{\alpha}(\cdot): $

$$
G_{\alpha}(x) = C_{11}(\alpha,d) \int_{R^d}
\frac{\exp(i x y)}{(1 + |y|^2)^{\alpha/2} } \ dy,
$$
which is used in the theory of Partial Differential Equations. \par
 As long as

$$
 \left| L_{\alpha} \ f \right|_p \le C_{12}(\alpha,d) |I_{\alpha} f|_p,
$$
see \cite{Adams1}, chapter 3, section 3.6, page 56,  we conclude that for the operator
$ L_{\alpha} $ are true  the {\it upper} estimations as for the operator $ I_{\alpha}. $ \par

Following,

$$
\left| L_{\alpha} \ f  \right|_q  \le
\frac{ C(\alpha,d) \ |f|_p }{[(p - 1)(d/\alpha - p)]^{1 - \alpha/d}};
$$
recall that $ 1/q = 1/p - \alpha/d, \ p \in (1,d/\alpha), \ q = q(p). $ \par
 Aside from we have from the Young inequality:

$$
\left| L_{\alpha} \ f  \right|_p \le C_1(\alpha,d) \ |f|_p,
$$
see \cite{Adams1}, p. 56. \par
It follows from the H\"older inequality that for all the values $ t $ from the
{\it closed} interval $ t \in [p, dp/(d - \alpha p)] $

$$
\left| L_{\alpha} \ f  \right|_t \le
\frac{ C(\alpha,d) \ |f|_p }{[(p - 1)(d/\alpha - p)]^{1 - \alpha/d}}.
$$
The last inequality may be rewritten as follows. We define the function
$$
\psi_0(t) = 1, \ t \in [p, dp/(d - \alpha p)], \ \psi(\cdot) \in
\Psi(p,  dp/(d - \alpha p));
$$
then

$$
|| L_{\alpha} \ f ||G(\psi_0) \le C_2(\alpha,d) \ \inf_{p \in (1, d/\alpha)}
\frac{C_1(\alpha,d) \ |f|_p}{ [(p-1)(d/\alpha - p)]^{1 - \alpha/d}}.
$$
 On the other hands, if for some $ p_0 \in (1, d/\alpha) \ |f|_{p_0} < \infty, $ then
 $ L_{\alpha} \ f \in G(\psi_0). $ \par

\vspace{3mm}

\vspace{3mm}

{\bf B.} Estimations of derivatives of Bessel's potential.\par

\vspace{3mm}

 We estimate in this subsection the norm in some Grand Lebesgue Space the
 derivatives of Bessel's potential function

 $$
 v^{(\xi)} = D^{\xi} \ v(x), \ v = L_{\alpha} \ f(x) = G_{\alpha}*f(x).
 $$
We suppose as before $ |\xi| < \alpha $ and define for the values
$ m \in (\max( 2\alpha/|\xi|,1), \infty)  $ the function

$$
\theta_{\alpha, \xi, d}(m) = \theta(m) = \frac{2 m |\xi|}{2m|\xi| - \alpha} \cdot
\psi^{|\xi|/2 \alpha} \left( \frac{2m|\xi|}{\alpha} \right)
\cdot \psi^{(\alpha - |\xi|)/(2 \alpha)}
\left(2m \left(1 -  \frac{|\xi|}{\alpha} \right) \right),
$$
where $ \psi(\cdot) \in \Psi(2|\xi|/\alpha, \infty). $ \par

\vspace{3mm}

{\bf Theorem 4.}

\vspace{3mm}

$$
||v^{(\xi)}||G(\theta) =
 ||D^{\xi}(G_{\alpha}*f ||G(\theta) \le C_3(\alpha, \xi, d) \ ||f||G(\psi). \eqno(28)
$$

\vspace{3mm}

{\bf Proof}. Let $ f \in G(\psi), $ otherwise is nothing to prove. We can and will
suppose $ |\xi| \ge 1 $ and $ ||f||G(\psi) = 1, $ then $ |f|_p \le \psi(p). $ \par
 We start from the inequality (see \cite{Adams1}, p. 57):

$$
|v^{(\xi)}(x) \le C \ (Mf(x))^{|\xi|/\alpha} \ |f(x)|^{1 - |\xi|/\alpha}. \eqno(29)
$$
We obtain after simple calculations using Cauchy - Schwartz inequality and denoting
for simplicity: $ A = \alpha \ p /(2 \ |\xi|): $

$$
\left[ (|v|_A)^A \right]^2 =
\left[\int_{R^d} \left|v^{(\xi)}(x) \right|^A \ dx \right]^2 \le \int_{R^d} |Mf(x)|^p \ dx \times
$$

$$
\int_{R^d} |f(x)|^{2A(1 - |\xi|/\alpha)} \ dx  \le C \ |f|_p^p \cdot
\left[ \frac{p}{p-1} \right]^p \times
$$

$$
\int_{R^d}|f(x)|^{p(\alpha - |\xi|)/|\xi|} \ dx  = C \ |f|_p^p \cdot
\left[\frac{p}{p-1} \right]^p \ \cdot
|f|^{p(\alpha - |\xi|)/|\xi|}_{p(\alpha - |\xi|)/|\xi|},
$$
therefore
$$
|v^{(\xi)}|_m =
\left|v^{(\xi)} \right|_A = \left|v^{(\xi)} \right|_{\alpha p/(2 |\xi|)} \le C \cdot
 \left( \frac{p}{p-1} \right)^{|\xi|/(2 \alpha)} \cdot
|f|_p^{|\xi/(2\alpha)} \cdot
|f|_{p(\alpha - |\xi|)/|\xi| }^{(\alpha - |\xi|)/(2 \alpha) } \le
$$

$$
C \cdot \frac{2m |\xi|}{2m |\xi| - \alpha} \cdot
 \psi^{|\xi|/(2\alpha)} \left( \frac{2m|\xi|}{\alpha} \right) \cdot
 \psi^{(\alpha - |\xi|)/(2 \alpha)} \left( 2 m \frac{\alpha - |\xi|}{\alpha }  \right)
 = C \ \theta(m).
$$
This competes the proof of the last theorem.\par

\bigskip

\section{The case of a bounded domain}

\vspace{3mm}

We consider in this section the  {\it truncated}  Riesz's operator

$$
u^{(B)} = u^{(B)}(x)  = I^{(B)}_{\alpha} f(x) = \int_B \frac{f(x - y) \ dy }{|y|^{d - \alpha} }, \eqno(30)
$$
where $ B $ is open bounded domain in $ R^d $ contained the origin and
 such that

$$
0 < \inf_{x \in \partial B } |x| \le \sup_{x \in \partial B } |x| < \infty, \eqno(31)
$$
$ \partial B $ denotes boundary of the set $ B. $ \par
It is known (see, e.g. \cite{Mitrinovich1}, p.90), that if
$ f \ \max(1, \log f) \in L_1(B), $ then $ f \in L_q(B). $ \par
 We can and will assume further without loss of generality that the set $ B $ is
 unit ball in the space $ R^d: $

 $$
B = \{x, \ x \in R^d, \ |x| < 1 \}. \eqno(32)
 $$
 Let us denote $ p^/ = p/(p - 1) $ and for the function $ \psi(\cdot) \in
 \Psi(1,d/\alpha) $ define

 $$
\nu(r) = \nu_{\psi}(r) = \inf_{p \in [1, d/(d - \alpha))}
\left[ \left( \frac{d }{d - \alpha } - p \right)^{-1 + \alpha/d} \cdot
  \psi \left( \frac{r p^/}{r + p^/} \right) \right]. \eqno(33)
 $$
 Note that the function $ \nu(r) $ may be simply estimated as follows: tacking
 the value $ p = p_0 = d/(d - \alpha) - C_1/r,  $ we conclude:

$$
\nu_{\psi}(r) \le C_2 \ r^{1 - \alpha/d} \ \psi( d/\alpha - C_3/r), \ r \ge C_4.
$$

\vspace{3mm}

{\bf Theorem 5}. Let $ \psi \in \Psi(1, d/\alpha). $  Then

\vspace{3mm}

$$
||I^{(B)}_{\alpha} f||G(\nu_{\psi}) \le C_6(\alpha,d) \ ||f||G(\psi). \eqno(34)
$$

\vspace{3mm}

{\bf Proof.}  As before, we suppose $ ||f||G(\psi) = 1,  $ hence

$$
|f|_p \le \psi(p), \ p \in [1, d/\alpha).
$$
We obtain by the direct computation for the values $ p $ from the
interval $ p \in [1, d/(d - \alpha) ) $:

$$
| \phi \cdot I(x \in B)|_p^p \asymp
\left( \frac{d }{d - \alpha } - p \right)^{-1},
$$
or equally

$$
| \phi \cdot I(x \in B)|_p \asymp
\left( \frac{d }{d - \alpha} - p \right)^{-1 + \alpha/d}.
$$
 As long as the function $ u^{(B)} $ may be written as a convolution

 $$
 u^{(B)} =  [\phi \cdot I(x \in B)]* f,
 $$
we can use the classical Young's inequality

$$
|u|_r \le | \ [\phi \cdot I(x \in B)] \ |_p \cdot |f|_k,
$$
where

$$
1 + \frac{1}{r} = \frac{1}{p} + \frac{1}{p},
$$
or equally
$$
 k = \frac{r p^/ }{r + p^/};
$$
and following

$$
|u|_r  \le C \ \left[ \left( \frac{d }{d - \alpha } - p \right)^{-1 + \alpha/d} \right]
 \ |f|_{rp^/ /(r + p^/)}. \eqno(35)
$$

We obtain substituting into (34) and using the equality $ ||f||G(\psi) = 1: $

 $$
|u|_r \le C \left[ \left( \frac{d }{d - \alpha } - p \right)^{-1 + \alpha/d} \cdot
  \psi \left( \frac{r p^/}{r + p^/} \right) \right]. \eqno(36)
 $$
Since the value $ p, \ p \in [1, d/(d - \alpha)) $ is arbitrary, we can minimize
the right - side of inequality  (35) over $ p: $

 $$
|u|_r \le C \inf_{p \in [1, d/(d - \alpha))} \
 \left[ \left( \frac{d }{d - \alpha } - p \right)^{-1 + \alpha/d} \cdot
  \psi \left( \frac{r p^/}{r + p^/} \right) \right] =
 $$

$$
C \nu_{\psi}(r) = C \nu_{\psi}(r) \ ||f||G(\psi).
$$
This completes the proof of theorem 5.\par
\hfill $\Box$

\vspace{3mm}

{\bf Example 3.} Let $ \psi \in G\Psi(1, d/\alpha; 0, \gamma); $ recall that this
means:

$$
\psi_{\gamma} (p) = (d/\alpha - p)^{-\gamma}, \ p \in [1, d/\alpha), \ \gamma =
\const > 0.
$$
Let  also $ f(\cdot) $ be a (measurable) function, $ f: R^d \to R $ such that

$$
|f|_p \asymp \psi_{\gamma}(p), \ p \in [1, d/\alpha).
$$
 This condition is satisfied, e.g., for the function

 $$
 f(x) = I(|x| > 1) \ \cdot \ |x|^{-\alpha} \ \cdot \ | \ \log |x| \ |^{\gamma - \alpha/d}, \ \gamma \ge \alpha /d.
 $$
 The  expression

$$
Z(p,r) =  \psi_{\gamma} \left( \frac{r p^/}{ ( r + p^/)} \right)
\cdot \left( \frac{d }{d - \alpha } - p \right)^{-1 + \alpha/d}
$$
achieves the asymptotical as $ r \to \infty $  minimal value at the point
$$
p_0 =  \frac{d}{d - \alpha} - 0.5 \left( \frac{d}{ d - \alpha } \right)^2 \
\frac{1}{r},
$$
and we obtain the following {\it upper} estimation for the $ L_r $ norm of the
function $ u(\cdot): $

$$
|u|_r \le C \ \inf_{p \in [1, d/(d - \alpha))} Z(p,r) \le C \ Z(p_0, r)
 \asymp r^{ 1 + \gamma - \alpha/d}, \ r \in [1, \infty). \eqno(37)
$$
Note that the last inequality (37) denotes that  the function $ u(\cdot) $
belongs to the  Orlicz's space over $ R^d $ equipped with the $ N \ - $ function

$$
N( t ) = \exp \left( |t|^m  \right) - 1, \  1/m = 1 + \gamma - \alpha/d, \ m \in
(0, \infty).
$$

\vspace{3mm}

 We intend now to establish the exactness of assertion of theorem 5,
 also in the one - dimensional case $ d = 1. $ \par

\vspace{3mm}

 {\bf Theorem 6.} The proposition (37)  is exact, i.e. there exists a function
 $ f_0 \in G(\psi_{\gamma}) $ for which

\vspace{3mm}

$$
|u|_r \ge C_2 \ r^{ 1 + \gamma - \alpha/d} = C_2 \ r^{1 + \gamma - \alpha }, \ r \ge 1.
$$

 Notice that the value $ r_0 = \infty $ in the unique "critical" point for the
 function $  r \to |u|_r, $ if $ f \in G(\psi_{\gamma}). $ \par
 {\bf Proof.} Let us denote $ \Delta = \gamma - \alpha/d.$ \par

As before, it is enough to investigate only the one - dimensional case $ d = 1. $ \par
 Let us choose
 $$
 f_0(x) = x^{-\alpha} \ |\log x|^{\Delta}, \ x \in (0,1/e)
 $$
and $ f_0(x) = 0 $ otherwise; $ \Delta = \const > 0.  $ \par
We find by the direct calculation: $ f_0 \in \cap_{p \in [1, 1/\alpha)} L_p $ and

$$
|f_0|_p \asymp (\alpha^{-1} - p)^{- \Delta - \alpha} =
(\alpha^{-1} - p)^{-\gamma},
$$
i.e. $ f_0 \in G(\psi_{\gamma}). $ \par

 Further, we have for the values $ x \in (0,1/e) $
$$
v_0(x) \stackrel{def}{=} I^{(B)}_{\alpha} f_0(x) =
\int_x^{1/e} |x - y|^{-\alpha} \ | \ \log|x - y| \ |^{\Delta} \ |y|^{\alpha - 1} \ dy
$$

$$
\sim \int_x^{1} |x - y|^{-\alpha} \ | \ \log|x - y| \ |^{\Delta} \ |y|^{\alpha - 1} \ dy
$$

$$
= \int_1^{1/x} \left[ |\log x| + |\log(z-1)| \right]^{\Delta} \ z^{\alpha - 1} \ dz \sim
\int_e^{1/x} \left[ |\log x| + |\log(z-1)| \right]^{\Delta} \ z^{\alpha - 1} \ dz
$$

$$
\sim |\log x|^{\Delta} \int_1^{1/x}(z-1)^{-\alpha} \ z^{\alpha - 1} \ dz \asymp
$$

$$
\sim |\log x|^{\Delta} \int_e^{1/x}(z-1)^{-\alpha} \ z^{\alpha - 1} \ dz \asymp
|\log x|^{\Delta} \int_e^{1/x} z^{-1} \ dz \sim |\log x|^{\Delta + 1};
$$

$$
|v_0|_r^r \asymp  \int_0^{1/e} |\log x|^{r(\Delta + 1)} \ dx  \sim
\int_0^{1} |\log x|^{r(\Delta + 1)} \ dx = \Gamma(r(\Delta + 1));
$$

$$
 \ |v_0|_r \sim  C_6 \ r^{\Delta + 1} = C_6 \ r^{1 + \gamma - \alpha}.
$$

This competes the proof of theorem 6. \par
\hfill $\Box$

\bigskip

\section{Concluding remarks}

\vspace{3mm}

{\bf A.} We consider in this subsection some generalization of the Riesz's potential operator of a view

$$
I_{\alpha, \beta} f(x) =
\int_{R^d} \frac{f(y) \  | \log |x - y| \ |^{\beta} \ dy }{ |x - y|^{d - \alpha} },
$$
$ \alpha = \const \in (0,d), \ \beta = \const > 0, $ or equally

$$
I_{\alpha, \beta} f(x) =
\int_{R^d} \frac{f(y) \  [ 1 + | \log |x - y| \ |]^{\beta} \ dy }{ |x - y|^{d - \alpha} },
$$

or more generally

 $$
 I_{\alpha, \beta}^{(S)} f(x) = \int_{R^d}
 \frac{f(y) \ | \log |x - y| \ |^{\beta} \ S(|\log|x - y| \ |) \ dy }
 { |x - y|^{d -\alpha} },
 $$
where $ \alpha = \const \in (0,d), \ \beta = \const > 0, $ and
$ S(z) $ is a {\it slowly varying} as $ z \to \infty $ continuous positive function:

$$
\forall \lambda > 0 \ \Rightarrow  \lim_{z \to \infty} S(\lambda z)/S(z) = 1.
$$

{\bf Lemma 2.}

$$
 |I_{\alpha,\beta} f|_q \le
 \frac{C \ |f|_p } {[(p - 1) \ (d/\alpha - p)]^{1 + \beta - \alpha/d } } \eqno(38)
 $$
{\it and the last inequality is sharp.}\par
 {\bf Proof } of the first assertion. We use the method described in
 \cite{Adams1}, pp. 49 - 54.\par
Let $ \chi =  \chi(z) $ be a positive at $ z \in (0,\infty) $ continuous decreasing
function such that $ \phi(\infty) \stackrel{def}{=} \lim_{z \to \infty} \phi(z) = 0. $
We define also a function

$$
\Phi(z) = \int_z^{\infty} \chi(t) \ dt,
$$
if there exists. \par
We have  for the values $ \delta \in (0, \infty) $ analogously to the assertion
in \cite{Adams1}, p. 49 - 51:

$$
\int_{y: |x - y| < \delta} \Phi(| x - y |) \ f(y) dy = \int_0^{\delta} \chi(r)
\int_{y: |x - y| < r} f(y) \ dy + \Phi(\delta) \int_{y: |x - y| < \delta} f(y) \ dy.
$$
 Without loss of generality we can assume that the function $ f(\cdot) $ is
 non - negative. \par
 As long as

 $$
 \int_{y: |x - y| \le r} f(y) \ dy \le C(d) \ r^d \ Mf(x),
 $$

 $$
 \int_{y: |x - y| \le \delta} f(y) \ dy \le C(d) \ \delta^d \ Mf(x),
 $$
we obtain the  estimate

$$
\int_{y: |x - y| < \delta} \Phi(|x - y |) \ f(y) dy \le C(d) \ Mf(x) \cdot
\left[\int_0^{\delta} r^d \ \chi(r) \ dr + \delta^d \ \Phi(\delta) \right]
$$

$$
\stackrel{def}{=}  Mf(x) \ A_{\chi, d}(\delta).
$$

Further, we have denote $ s = p/(p-1) $ and use the H\"older inequality:

$$
\int_{y: |x - y| \ge \delta} \Phi(|x - y|) \ f(y) \ dy \le |f|_p \cdot
\left[ \int_{y: |x - y| > \delta} \Phi^s(|x - y| \ dy  \right]^{1/s} =
$$

$$
C_1(d) \ |f|_p \left[\int_{\delta}^{\infty} r^{d-1} \Phi^s(r) \ dr \right]^{1/s}
= D_{\chi}(p, \delta) \ |f|_p,
$$
where

$$
D_{\chi}(p, \delta) \stackrel{def}{=} C_1(d) \left[\int_{\delta}^{\infty} r^{d - 1} \
\Phi^s(r)  \ dr \right]^{1/s},
$$
if there exists for some values $ s $ from some non - trivial interval
$ s \in (d/(d - \alpha), s_0); $ if $ s_0 < \infty, $ then we define formally
$ D_{\chi}(p, \delta) = + \infty.  $\par

We conclude tacking into account the partition

$$
\int_{R^d} \Phi(|x - y|) \ f(y) \ dy = \int_{y: |x - y| < \delta} \Phi(|x - y|) \ f(y) \ dy
$$

$$
+ \int_{y: |x - y| \ge \delta} \Phi(|x - y|) \ f(y) \ dy:
$$

$$
\int_{R^d} \Phi(|x - y|) \ f(y) \ dy \le Mf(x) \ A_{\chi, d}(\delta) +
D_{\chi}(p, \delta) \ |f|_p.
$$
Therefore,

$$
\int_{R^d} \Phi(|x - y|) \ f(y) \ dy \le \inf_{\delta > 0}
\left[ Mf(x) \ A_{\chi, d}(\delta) +D_{\chi}(p, \delta) \ |f|_p \right]
$$

$$
\stackrel{def}{=} H(p, Mf(x), |f|_p).
$$
 Solving the last inequality, we obtain denoting
 $$
 w(x) =  \int_{R^d} \Phi(|x - y|) \ f(y) \ dy:
 $$
the inequality of a view

 $$
 G(p, w(x), |f|_p) \le (Mf(x))^p,
 $$
 and after the integration

$$
 \int_{R^d} G(p, w(x), |f|_p) \ dx \le |Mf(x)|_p^p \le C^p(\alpha,d) \ |f|_p^p \ (p-1)^{-p},
 $$

 $$
  \left[\int_{R^d} G(p, w(x), |f|_p) \ dx \right]^{1/p}
   \le |Mf(x)|_p \le C(\alpha,d, \chi) \ |f|_p \ (p-1).
 $$

Since the relation between  the functions $ f(\cdot) $ and $ w(\cdot) $ is linear,
the last inequality has a view

$$
|w|_q \le C_2(\alpha,p, \chi(\cdot)) \ |f|_p.
$$

 Choosing the function $ \Phi(r) $ as  follows:

 $$
 \Phi(r)= r^{\alpha - d} \ | \log r|^{\beta},
 $$
 we obtain the first assertion of lemma 2.\par

\vspace{3mm}

{\bf Proof} of the second assertion of the lemma 2. \par
 Notice that it may be proved  a more general assertion:

 $$
 |I_{\alpha, \beta}^{(S)} f|_q \le
 \frac{C \ |f|_p }{[(p - 1) \cdot (d/\alpha - p)]^{1 + \beta - \alpha/d} }
  \cdot S( q(d - \alpha) - d)^{-1}) \ S(q),
 $$
and the last estimation is asymptotically exact: at the same examples as at the
proof theorem 2 can be used by the proof of inverse inequalities. \par
  More detail: let $ d = 1, \ p \in (1, 0.5(1 + d/\alpha)), \ p \to 1 + 0, $
 or equally $ q \to 1/(1 - \alpha) + 0 $ and let us choose

 $$
 f(x) = x^{-1} \ (\log x)^{\Delta} \ I(x > e), \ \Delta = \const \ge 0.
 $$
The value $ |f|_p $ was calculated before, let us estimate the $ L_q $ norm
 of a function
 $$
 u(x) =  I_{\alpha, \beta}^{(S)} f(x).
 $$
 We have using the properties of slowly varying functions
  as $ x \to \infty, \ x > e: u(x) = $

$$
\int_e^{\infty} |x-y|^{\alpha-1} \ y^{-1} \ (\log y)^{\Delta} \ |\log|x-y| \ |^{\beta} \
S(|\log|x-y||) \ dy = x^{\alpha - 1} \times
$$

$$
\int_{e/x}^{\infty} z^{-1} \ |z-1|^{1 - \alpha} \ x^{\alpha - 1} \
[\log x + \log z]^{\Delta} \
[\log x + \log(|z-1|)]^{\beta} S( \log x + \log(|z-1|)) \ dz   \asymp
$$

$$
x^{\alpha - 1} \ [\log x]^{\Delta + \beta} \ S(\log x) \int_{e/x}^{\infty}
z^{-1} \ |z - 1|^{\alpha - 1} \ dz \asymp
$$

$$
x^{\alpha - 1} \ [\log x]^{\Delta + \beta} \ S(\log x) \int_{e/x}^{1/e}
z^{-1} \ |z - 1|^{\alpha - 1} \ dz \asymp
$$

$$
x^{\alpha - 1} \ [\log x]^{\Delta + \beta} \ S(\log x) \int_{e/x}^{1/e}
z^{-1} \ dz \asymp x^{\alpha - 1} \ [\log x]^{\Delta + \beta + 1} \ S(\log x);
$$
following,

$$
|u|_q^q  \asymp \int_e^{\infty} x^{-q(1 - \alpha)} \ |\log x|^{q(\Delta + \beta+1)} \
S^q(|\log x|) dx  \asymp
$$

$$
\int_0^{\infty} e^{-y[q(1-\alpha) - 1]} \ y^{q(\Delta + \beta + 1)} \ S^q(y) \ dy =
$$

$$
[q - 1/(1 - \alpha)]^{-q(\Delta + \beta + 1) - 1} \ \int_0^{\infty} e^{-z} \
 z^{q(\Delta + \beta + 1)} \ S^q \left( \frac{z}{q(1 - \alpha) - 1  }\right) \ dz =
$$

$$
[q - 1/(1 - \alpha)]^{-q(\Delta + \beta + 1) - 1} \
S^q \left( \frac{ 1 }{q - 1/(1 - \alpha)}  \right).
$$
Therefore  we have as $ q \to 1/(1 - \alpha) + 0 $

$$
|u|_q \asymp [q - 1/(1 - \alpha)]^{-\Delta - \beta - 2 + \alpha} \
S \left( \frac{1 - \alpha}{q(1 - \alpha) - 1} \right) \asymp
$$

$$
 [p - 1]^{-\Delta - \beta - 2 + \alpha} \
S \left( \frac{1 - \alpha}{q(1 - \alpha) - 1} \right).
$$
 Substituting into the expression for the left hand side of inequality (48), we
 can see that the relation between the left hand side to the right hand side
 calculated for the function $ f $ is bounded from below as $ p \to 1 + 0. $ \par
  The case $ p \to 1/\alpha - 0 $ or equally $ q \to \infty $ may be considered
  analogously, by means of example of a function

  $$
  g(x) = x^{-\alpha} \  |\log x|^{\Delta} \ I(x \in (0, 1/e)).
  $$
We have in this case:

$$
|g|_p \asymp (1/\alpha - p)^{-\Delta - \alpha};
$$

$$
v(x):=I_{\alpha,\beta}^{(S)} g(x) \asymp |\log x|^{\Delta + \beta + 1}
\ S(|\log x|), \ x  \to 0+;
$$

$$
|v|_q \asymp q^{\Delta + \beta + 1} \ S(q) \asymp (1/\alpha - p)^{-\Delta - \beta - 1}
 \ S \left( \frac{1}{1/\alpha - p }  \right)
$$
etc. \par

\vspace{4mm}

 Let us denote for arbitrary function $ \psi(\cdot) $ from the class
 $ \Psi(1, d/\alpha) $

  $$
  \zeta_{\alpha,\beta}^{(S)}(q) =
  \frac{\psi(p) \ S((p-1)^{-1}) \ S( q(d - \alpha) - d)^{-1}}
  { [(p - 1) \ (d/\alpha - p) ]^{1 + \beta - \alpha/d }}.
  $$
Here as before
$$
p \in (1, d/\alpha), \ q^{-1} = q^{-1}(p) = (p^{-1} - \alpha/d)^{-1}
\in (d/(d - \alpha), \infty).
$$
 We find  analogously the assertion of theorem 1: \par

 \vspace{3mm}

 {\bf Theorem 7.}

\vspace{3mm}

 $$
 ||I_{\alpha,\beta}^{(S)} f ||G \left(\zeta_{\alpha,\beta}^{(S)} \right) \le
 C \left(\alpha,\beta, d, S(\cdot) \right) \ ||f||G(\psi).
 $$

\vspace{4mm}

{\bf B.} Let us consider in this subsection the  {\it generalized truncated}  Riesz's operator

$$
u^{(B)}_{\beta} = u^{(B)}_{\beta}(x)  = I^{(B)}_{\alpha, \beta} f(x) =
\int_B \frac{f(x - y) \ |\ \log |y| \ |^{\beta} \ dy }{|y|^{d - \alpha} },
$$
where $ \beta = \const > 0, $
$ B $ is open bounded domain in $ R^d $ contained the origin and
 such that

$$
0 < \inf_{x \in \partial B } |x| \le \sup_{x \in \partial B } |x| < \infty,
$$
$ \partial B $ denotes boundary of the set $ B. $ \par
 We can and will assume further without loss of generality that the set $ B $ is
the ball of a radius $ 1/e $ in the space $ R^d: $

 $$
B = \{x, \ x \in R^d, \ |x| < 1/e \}.
 $$
 Let us denote for the function $ \psi(\cdot) \in
 \Psi(1,d/\alpha) $ define

 $$
\nu^{(\beta)}(r) = \nu_{\psi}^{(\beta)}(r) = \inf_{p \in [1, d/(d - \alpha))}
\left[ \left( \frac{d }{d - \alpha } - p \right)^{-1 - \beta + \alpha/d} \cdot
  \psi \left( \frac{r p^/}{r + p^/} \right) \right]. \eqno(39)
 $$

\vspace{3mm}

{\bf Theorem 8}. Let $ \psi \in \Psi(1, d/\alpha). $  Then

\vspace{3mm}

$$
||I^{(B)}_{\alpha, \ \beta} f||G(\nu_{\psi}^{(\beta)}) \le
C_9(\alpha,d) \ ||f||G(\psi). \eqno(40)
$$
The {\bf proof} is at the same as in theorem 5; it based on the following
equality on the function

$$
\phi_{\beta}(x) = I(x \in B) \ |x|^{\alpha - d} \ | \ \log |x| \ |^{\beta}:
$$

$$
| \ \phi_{\beta}(\cdot) \ |_p \asymp
\left[ \frac{d}{d - \alpha} - p \right]^{- \beta - 1 + \alpha/d }, \eqno(41)
$$
$ p \in [1, d/(d - \alpha));$ hence when $ \psi \in \Psi(1, d/\alpha), $ then

 $$
\nu^{(\beta)}(r) = \nu_{\psi}^{(\beta)}(r) = \inf_{p \in [1, d/(d - \alpha))}
\left[ \left( \frac{d }{d - \alpha } - p \right)^{-\beta - 1 + \alpha/d} \cdot
  \psi \left( \frac{r p^/}{r + p^/} \right) \right].
 $$

 As in the section 3, we conclude that the estimation (40) is sharp. \par
We can consider the more general case when the integral operator has a view:

$$
I^{(B, S)} f = \phi^{(S)}_{\alpha, \beta}*f, \ \alpha \in (0,d), \ \beta \ge 0,
 \eqno(42)
$$
where

$$
\phi^{(S)}_{\alpha, \beta}(x) = I(x \in B) \cdot |x|^{\alpha - d} \cdot
| \log |x| \ |^{\beta} \ S(| \log |x| \ |), \eqno(43)
$$
$ S(z) $ is a slowly varying as $ z \to \infty $ continuous positive function.\par

Let us denote

$$
\nu^{(S)}_{\alpha, \beta}(r) = \inf_{p \in [1, d/(d - \alpha)}
\left[ \left( \frac{d }{d - \alpha } - p  \right)^{-1 - \beta + \alpha/d} \cdot
S \left( \frac{d - \alpha }{d - p(d - \alpha) }  \right) \cdot
\psi \left( \frac{r p^/}{r + p^/ }   \right)     \right].
$$

We assert analogously to the theorem 7: \\

{\bf Theorem 9.}
$$
||I^{(B,S)}_{\alpha, \ \beta} f||G(\nu^{(S)}_{\alpha,\beta}) \le
C_9(\alpha,d) \ ||f||G(\psi). \eqno(44)
$$
{\bf Proof.} It is sufficient for the proof of the last assertion to calculate all the moments of a function

 $$
 R(x) = I(x \in B) \ |x|^{-\alpha}  \ | \ \log |x| \ |^{\Delta} \ S(|\ \log |x| \ |).
 $$
We have for the values $ p $ from the set $ p \in [1, d/\alpha) $
(for the values $ p $ greatest than $ d/\alpha: \ p > d/\alpha  \ |R|_p = \infty ) $
using the multidimensional spherical coordinates:

$$
|R|_p^p \asymp \int_0^1 r^{d - 1 - \alpha p} \ |\log r|^{\Delta p} \
S^{p}(|\log r|) \ dr =
\int_0^{\infty} e^{-y(d - \alpha p)} \ y^{\Delta p} \ S^p(y) \ d y =
$$

$$
(d - \alpha p)^{-\Delta p - 1} \ \int_0^{\infty} e^{-z} z^{\Delta p} \
S^p \left( \frac{ z }{ d - \alpha p } \right) \ dz \sim
(d - \alpha p)^{-\Delta p - 1} \times
$$

$$
 \ S^p \left( \frac{ 1 }{ d - \alpha p } \right) \
\int_0^{\infty} e^{-z} \ z^{\Delta p} \ dz =
(d - \alpha p)^{-\Delta p - 1} \ S^p \left( \frac{ \alpha }{ d/\alpha - p } \right) \
\Gamma (\Delta p + 1),
$$
as long as the function $ S(\cdot) $  is slowly  varying. \par
Further,

$$
|R|_p \asymp (d - \alpha p)^{-\Delta  - \alpha/d} \
S \left( \frac{ \alpha }{ d/\alpha - p } \right).
$$

 \vspace{3mm}

{\bf B.}\par
We consider now  a so - called fractional sublinear  maximal operator:
$$
M_{\alpha}f = \sup_{x \in R^d} \ \sup_{\rho > 0} \left[ \rho^{\alpha - d}
\int_{y: |x - y| \le \rho} |f(y)| \ dy \right].
$$
As long as

$$
C_{10}(\alpha,d) \ | M_{\alpha} f|_p \le |I_{\alpha} f|_p \le
C_{11} (\alpha,d) \ | M_{\alpha} f|_p,  \ p \in (1, d/\alpha),
$$
see \cite{Adams1}, chapter 3, section 3.6, we conclude that for the operator
$ M_{\alpha} $ are true at the same estimations as for the operator $ I_{\alpha}. $ \par

\vspace{3mm}

\vspace{3mm}

{\bf E.} The $ L_p \to L_q $ estimations  for many examples of integral operators
(regular and singular) in the {\it weighted } Lebesgue spaces $ L_p(R^d,w) $ with
exact values of its norm see, for example, in \cite{Perez1}. \par
 This estimations allow to generalize described results on the weighted spaces.\par

\end{document}